\documentclass[10pt]{amsart}
\usepackage[cp1251]{inputenc}
\usepackage[english]{babel}
\usepackage{amsmath}
\usepackage{amssymb}
\usepackage{amsfonts}
\usepackage{amscd}

\newcommand{\ii}{\infty}
\newcommand{\ee}{\varepsilon}
\newcommand{\aaa}{\alpha}
\renewcommand{\lll}{\langle}

\newcommand{\rrr}{\rangle}

\begin{document}

\title{Obstructions for uniform stability of $C_0$-semigroup}

\author{K. V. Storozhuk}

\address{Sobolev Institute of mathematics SB RAS}
\begin{abstract}
Let  $T_t:X\to X$ be a $C_0$-semigroup with generator $A$.  We prove
that if the abscissa of uniform boundedness of the resolvent
$s_0(A)\geq 0$ then for each a non-decreasing function $h(s):\Bbb
R_+\to R_+$, there are $x'\in X'$ and $x\in X$ such that $\int_0^\ii
h(|\lll x',T_tx\rrr|)dt=\ii$. If $i\Bbb R\cap Sp(A)\neq \emptyset$
then such $x$ may be taken in $D(A^\infty)$.
\end{abstract}

%\begin{keyword}
%exponentially stable operator semigroup (keyword)

%\end{keyword}

\maketitle

\section{Introduction and preliminary results}

Let $X$ be a Banach space. $C_0$-semigroup  $T_t:X\to X$
 is called uniformly exponentially stable  (UES), if $\|T_tx\|$ decays
exponentially with $t$ for all $x\in X$. By the Uniform Boundedness
Principle (UBP) it is equivalent to  $\|T_t\|\to_{t\to\ii} 0$.

In the finite-dimensional case, these conditions are equivalent to
the fact that $\|T_tx\|\to 0$ as $t\to\ii$ for all $x\in X$. The
basic infinite-dimensional counterexample is given by the semigroup
of shifts on $L_2(\Bbb R_+)$. Here $\|T_t\|\equiv 1$ for all $t$,
but $\|T_tx\|\to 0$ for all~$x$. However, the absence of (UES) for a
semigroup implies the existence of the vectors whose orbits, while
tending to zero, do it ``very slowly''. In this connection, we
mention the articles \cite{1Mu, 2N3} which, in particular, imply
that if the spectral radius $r(T)=1$, then, for every sequence
$1>\aaa_n\to 0$, there exists $x\in X$ such that $\|T^nx\|>\aaa_n$
for all $n$.

 For $C_0$-semigroups, let us give a typical result demonstrating
the absence of the upper estimate for the decrease rate "in the
integral sense".

\smallskip

\bf Theorem 0. \it Suppose that $\forall t\geq 0\ \|T_t\|\geq 1$.
For any non-decreasing function $h:\Bbb R_+\to R_+$ (``\it good\rm''
function), there exists $x\in X$ such that $\int_0^\ii h(\|T_t
x\|)dt=\ii.$ Moreover, if the semigroup $T_t$ is unbounded then for
some $x\in X$ $\|T_tx\|>1$ for  $t$ from a set of infinite measure.

\smallskip

Analogous results hold also for a more general situation of
evolution operators $U(t,s)U(s,r)=U(t,r):X\to X$. There are several
proofs of such results going back to \cite{3Da, 4Pa, 5Ro}.
%
%The use of the formula (**)  is a quite analogous to beginning of
%the proof of the Theorem 3.2.2 in \cite{N1}.
%
In a short proof the first part belongs to the author \cite{6St},
the second part is  quite analogous to the beginning of the proof of
the Theorem 3.2.2 in \cite{7N1}. We will present this proof, first,
because of its brevity, and, second, because we will compare other
reasoning with it.
\smallskip

 \it Proof of Theorem 0\rm. The following ``backward''
 estimate holds for a bounded $C_0$-semigroup:
$$\text{Let }\ C=\sup_{t>0}\{\|T_t\|\}.\ \text{Then }\ \forall t_0\ \forall t\in[0,t_0],\
 \forall x\in X\ \|T_tx\|\geq \frac{\|T_{t_0}x\|}{C}\eqno(*) $$

 Choose a sequence $\aaa_n\to 0$ that decreases so slowly that
$n h(\aaa_n)\to\infty$. We have $\lim\frac{\|T_n\|}{\aaa_n}=\ii$.
Following the UBP, choose $x\in X$ such that
$\overline{\lim\limits_{n\to\infty}}\frac{\|T_nx\|}{\aaa_n}\geq C$.
Then
$$\int\limits_0^\ii h(\|T_t x\|)dt=\sup_{n} \int\limits_0^{n} h(\|T_t
x\|)dt\overset{(*)}{\geq} \sup_n\int\limits_0^{n}
h\left(\frac{\|T_{n} x\|}{C}\right)dt \geq \sup_n n\cdot
h(\aaa_{n})=\infty.
$$

Now, let  $\|T_t\|$ be an unbounded $C_0$-semigroup. In this case,
we have a weaker analog of $(*)$, a ``finite-time backward
estimate''. For example,
$$\text{Let }\ C=\sup_{t\in[0,1]}\{\|T_t\|\}.\ \text{Then }\ \forall t_0\ \forall t\in[t_0-1,t_0],\
\forall x\in X\ \|T_tx\|\geq \frac{\|T_{t_0}x\|}{C}\eqno(**) $$ The
estimate $(**)$ is weaker than the estimate ($*$). On the other
hand, according to UBP, for the unbounded semigroup there exists
$x\in X$ whose orbit $T_tx$ does not tend to zero. Let
$\|T_{n_k}x\|>C$, $n_k\to\infty$. Formula $(**)$ implies that
$\|T_tx\|>1$ with $t$ in the set
 $\cup_{k=1}^\ii[n_k-1,n_k]$.
\smallskip

The questions, concerning the orbits' slow approaching to zero in
the weak topology of the space  $X$, appeared for the first time and
started to be discussed in  \cite{8PZ,9Hu, 10We}. A review of this
subject can be found in \cite{7N1}. A question, which is analogous
to the conclusion of Theorem 0 for the weak topology, is the
following: under which assumptions on the semigroup one can claim
that for each "good" function $h$
$$ \exists x\in X,\ x'\in X' \ \int_0^\ii h(|\lll x',
T_tx\rrr|)\,dt=\ii. \eqno(1) $$

The absence of UES alone is not sufficient here, as can be shown by
the example (see  \cite{11GVW}; \cite{12NSW}, example~1.5) of the
semigroup of shifts on $L^1(\Bbb R_+, e^{t}dt) \cap L^p(\Bbb R_+).$
This semigroup is not UES but is weakly $L^1$-stable; i.e., $$
\forall x\in X,\ x'\in X' \ \int_0^\ii|\lll x', T_tx\rrr|\,dt<\ii.$$

From the standpoint of geometry, this example looks rather
 surprising: some orbits are "far" from zero; at the same time,
 the orbit of each vector $x$ spends almost all the time
 "arbitrarily close" to each hyperplane $(\ker x')$!

It is known that the conclusion (1) holds, for example, for the
bounded $C_0$-semigropus, if one requires the absence of not only
the \it uniform exponential stability \rm but simply the \it
exponential stability\rm.
 (see the beginning of the next section and Proposition 1).

The main result of the present article (Theorem 1) is the estimates
from which it follows, in particular, that, in the Proposition 1,
the condition of the semigroup being bounded is dispensable.

In the next section we briefly discuss some asymptotic parameters of
the semigroup which are necessary for precise formulation of the
results and formulate the appropriate statements. The third section
is a proof of the basic result. In the fourth section we prove that,
under some natural assumptions, the vector $x$ in the formula (1)
can be chosen among "smooth" vectors.

\section{Asymptotic parameters and the basic
results}\smallskip%

Uniform growth bound of semigroup is
$$\omega_0(T)=\lim_{t\to\ii}\frac{\ln\|T_t\|}{t}=\sup\{\overline{\lim}_{t\to\ii}\frac{\ln\|T_tx\|}{t}
\mid x\in X\}.$$

Let $A$ be a generator of the semigroup and let $D(A)$ be a domain
of $A$. The growth bound of semigroup is
$$\omega_1(T):=\sup\{\overline{\lim}_{t\to\ii}\frac{\ln\|T_tx\|}{t}
\mid x\in D(A)\}.$$

Let us explain informally the meaning of $\omega_0$ and $\omega_1$.
Consider an abstract Cauchy problem $\frac{dz}{dt}=Az,\ z(0)=x\in
X$. The function $z(t)=T_tx$ is called a \it mild solution \rm of
ACP. If the initial value $x$ lies in $D(A)$, then it is natural for
the solution to be called \it classical\rm, or \it smooth \rm
solution (accordingly, such initial values $x$ from $D(A)$ are
called the smooth vectors). So, the parameters $\omega_0$ and
$\omega_1$ describe the growth of the mild (accordingly, smooth)
solutions of our ACP.

Clearly, $\omega_1\leq\omega_0$. It is easy to see that a semigroup
is UES if and only if $\omega_0(T)<0$. If $\omega_1(T)<0$, then the
semigroup is simply called exponentially stable, or (ES)

For the semigroup of shifts, mentioned in the end of the
introduction, $\omega_1<0=\omega_0$; therefore though it is not UES
but, nevertheless, it is ES.

Theorem 0 can still be "rescued" for the weak topology also, if, for
example, one requires that the semigroup, together with being not
UES, is not ES.
\smallskip

\bf Proposition 1 \it  (see \cite{7N1}, Theorems~4.6.3(i) and
~4.6.4). If a uniformly bounded $C_0$-semigroup $T_t:X\to X$ is not
expo\-nen\-ti\-ally stable, then the following holds:

(1) for any ``good'' $h>0$, there is $x\in X $ such that $\int_0^\ii
h(|\lll x',T_tx\rrr|)dt=\ii; $

(2) there exists an $\ee> 0$ such that for each $m > 0$ there exist
norm-one vectors $x\in X$ and $x'\in X'$ such that $m<mes\{t\mid
|\lll x',T(t)x\rrr|\geq\ee\}.$\rm

Let us describe two more spectral characteristics, between which the
growth $\omega_1$  of the semigroup is "confined". It will be
necessary for us to formulate of basic results of article and
strengthen Proposition 1. It is the spectral bound $s(A)$ and the
abscissa $s_0(A)$ of uniform boundedness of the resolvent of $A$. We
have the following diagram of inequalities
 $\begin{array}{ccc}
      \omega_1 &\!\!\!\!\! \leftarrow  & \!\!\!\! \omega_0 \\
    \downarrow &\!\!\!\!\! \nwarrow    & \!\!\!\! \downarrow\\
             s &\!\!\!\!\! \leftarrow  & \!\!\!\!  s_0 \\
  \end{array}$ (here ``$\leftarrow$'' stands for ``$\leq $'').

In Hilbert spaces, $s_0=\omega_0$ (\cite{13Ge}). For positive
semigroups, $s=\omega_1=s_0$; moreover, $s$ is reached by
$\sigma(A)$, i.e., there exists $\lambda\in \sigma(A),\
re(\lambda)=s$ (\cite{11GVW,14Na,15Nb}).

On the other hand, for each arrow, there is an example of strict
inequality. The first historically  example of Foia\c{s}'
\cite{16Fo}, in which $-2=s<s_0=\omega_0=0$ is known rather badly.
In the foundation of the majority of other such examples there lies
an example due to Zabczyk \cite{17Zb}. (Note that Zabczyk himself
refers to Foya\c{s} in the text of the paper \cite{17Zb}). In
Zabczyk's example $s<\omega_1=\omega_0$; $s$ is reached and
$\|T_t\|=e^{t\omega_0}$. These and other examples can be found in
\cite{14Na} and \cite{7N1}.
\smallskip

It turns out that the conclusion of Proposition 1 is valid also for
non-bounded semigroups, whereas the condition of "being not ES"
($\omega_1(T)\geq 0$) can be substituted by a weaker condition
"$s_0(A)\geq 0$". This follows from our main result:
\pagebreak%\bigskip\rm

\bf Theorem 1\it.
  Let $A$ be the generator of $C_0$-semigroup
$T_t:X\to X$ and let $s_0(A)\geq 0$. For any two sequences $0<m_1<
m_2<\cdots$ and $\gamma_k>0$, $\gamma_k\to 0$ there is $x'\in X'$
$x\in X$ and a family of sets $U_k\subset \Bbb R_+$ such that
$$\forall k\in \Bbb N\ \mu(U_k)\geq m_k,\
 \forall t\in U_k\ \ |\lll x',T_tx\rrr|>\gamma_k. %$ \frac{5}{10^{2^k-1}}.
$$\rm

To derive Proposition 1 from Theorem 1, it suffices to find, for the
function $h$, the numbers $m_k$ growing so rapidly that $m_k\cdot
h(\gamma_k)\to\infty$,  in which case the integral in Proposition 1
diverges (compare with the proof of Theorem 0).

\it Notice\rm. In fact, the proof of Theorem~4.6.3 of \cite{7N1}
uses only the condition on $s_0$ instead of $\omega_1$. However, the
assumption of boundedness is essential for those methods (the proof
in \cite{7N1} is based on the technique of rearrangement-invariant
Banach function spaces). In our proof of Theotem 1 we will need some
methodological and technical tricks which would be unnecessary under
the assumption of the boundedness of the semigroup.

 \section{Proof of Theorem  1}

Our proof of Theorem 1 is based on the following lemma. \smallskip

 \bf
Lemma 1\it. Suppose that $s_0(A)=0$. Then, for all $\delta>0$ and
$t_0<\ii$ there are $\beta=\beta(\delta, t_0)\in \Bbb R$ and $y\in
D(A)$ such that $\|y\|=1$, $\|Ay\|\approx |\beta| $ and
$\|T_ty-e^{i\beta t}y\|<\delta$ for all $t\in [0,t_0]$. Moreover,
such  $y$ may be chosen in $D(A^\ii)=\cap_{n}D(A^n)$.\rm
\smallskip

\it Notice\rm. Geometrically, Lemma 1 means that there is a unit
vector $y$ staying near the (complex) line $\Bbb{C}y$ for a long
time (at that, this vector is "spinning" in the corresponding real
plane with angular velocity  $\beta$ ).
 At the same time, the assertion~(2) of Proposition 1 means merely that there is a unit vector $x$
staying away from the hyperplane $\ker\, x'$ for a long time. So,
(2)  (proved in \cite{7N1} with the use of assertion~(1)) follows
already from Lemma~1.

\it Proof of Lemma 1\rm.
 Note that, given a semigroup $T_t$
with generator $A$, we have
$$\forall x\in D(A) \ \forall t>0, \beta\in\Bbb R\ \
\|T_tx-e^{i\beta t}x\|\leq t\cdot\sup_{s\in[0,t]}\|T_s\|\cdot
\|(A-i\beta) x\|.\eqno(2)$$ Indeed, $(A-i\beta)$ is the generator of
the semigroup $e^{-i\beta t}T_t$, so $\|T_tx-e^{i\beta
t}x\|=\|e^{-i\beta t}T_tx-x\|=\|\int_0^tT_s(A-i\beta)x\,ds\|$ etc.

As $s_0(A)=0$, the resolvent $(A-\lambda)^{-1}$ of operator $A$ is
unbounded if $\lambda$ is near the imaginary axis. Choose
$\lambda_n=\alpha_n+i\beta_n\in \Bbb C$, $\alpha_n\to 0$,
$\|(\lambda_n-A)^{-1}\|\to\infty$. Take the vectors $y_n\in D(A)$
such that $\|y_n\|=1$ but $(A-\lambda_n)y_n\to 0$. In particular,
$(A-i\beta_n)y_n\to 0$. Take $n$ so large that
$\|(A-i\beta_n)y_n\|<\frac{\delta}{t_0\cdot\sup_{s\in[0,t]}\|T_s\|}$.

The set $D(A^\ii)$ is dense in $D(A)$ in the graph norm
(\cite{14Na}, 1.9(iii)), so such $y$ may be chosen in $D(A^\ii)$.
 Now,
involving (2), we finish the proof of the lemma.$\square$\smallskip

\it Proof of Theorem 1\rm. Clearly, if there are $x'\in X'$,  $x\in
X$, and $\delta>0$ such that $mes\{t>0\mid |\lll x',T_tx
\rrr|>\delta\}=\infty$, then there is nothing to prove. In
particular, we may assume from the very beginning that
$$\forall x\in D(A)\ \forall x'\in X'\ \forall \delta>0 \ \  mes\{t>0\mid |\lll x',T_tx
\rrr|>\delta\}<\ii. \eqno (3)$$

Suppose that $X$ is separable, this will enable us to use the
sequential compactness of  $X'$  in the $*$-weak topology.  It does
not restrict the generality: it is not difficult to see that in
general case Theorem 1 can be applied to an appropriate separable
subspace $X_s\subset X$,  invariant under the action of the
semigroup, and then continue the corresponding functional $x_s'\in
(X_s)'$  to a functional $x'$ on the entire $X$.

Last, note that it is enough for us to prove the theorem for any
concrete sequence $\gamma_k'\to 0$. (Let us explain why it does not
restrict the generality either. Let $m_k\in \Bbb N$ and $\gamma_k\to
0$. We may suppose that $\gamma_k<\gamma_1'$ for all $k$. Put
$n(k)=\max\{n\mid \gamma_n\geq\gamma_k'\}$ and
$m_k'=m_1+\cdots+m_{n(k)}$. It is easy to see that if $x, x'$
satisfy the hypothesis of Theorem 1 with $m_k'$ and $\gamma_k'$ then
they will also satisfy of Theorem 1 for initial
$m_k,\gamma_k$.)\smallskip

We prove Theorem 1 for $\gamma_k=\frac{5}{10^{2^k-1}}$.

Following Lemma 1, choose a sequence of numbers $\beta_n\in\Bbb R$
and vectors $ y_n\in D(A)$ such that $\| y_n\|=1$ and $\|T_t
y_n-e^{i\beta_nt} y_n\|\leq{\frac{1}{10}}$ for some $\beta_n\in\Bbb
R$ and all $t\in [0,n]$.

For each $y_n$, choose any dual $ y'_n\in X'$, $\| y'_n\|=\lll
 y'_n,  y_n\rrr=1$. The sequence $ y'_n$ contains a subsequence
that *-weakly converges to some $ y'\in X'$. We may consider $
y_n'\overset{\sigma^*}{\to}  y'$. It easy to see that
$$\forall t\in[0,n]\ \frac{9}{10}<\|T_t y_n\|<\frac{11}{10}, \  \frac{9}{10}<|\lll  y_n',T_t y_n\rrr|<\frac{11}{10}
\eqno(4)$$
 We shall construct $x$ and $x'$ as the limits of $x_k$ and $x'_k$,
$x_k=\sum_{l=1}^k \frac{\pm y_{n_{l}}}{10^{2^{l-1}-1}}, \
x_k'=\sum_{l=1}^k \frac{\pm y'_{n_{l}}}{10^{2^{l-1}-1}},$ where the
numbers $n_{l}$ and the signs $\pm$ are to be found.\smallskip

{\bf The construction of the vector $x_1$}. Let $n_1\in \Bbb N$,
$n_1\geq m_1$. Put $U_1=[0,n_1]$. Put $x_1':= y_{n_1}',$ $x_1:=
y_{n_1}$. Then
$$\forall t\in U_1\ |\lll
x_1',T_tx_1\rrr|>\frac{9}{10}\eqno(5)$$

{\bf The construction of the vector $x_2$}.

 The assumption, expressed by the formula ~(3),  allows choosing a
 sufficiently large but compact set $ \tilde {U} _ {2} $ such that
$\mu(\tilde{U}_{2})\geq 4m_2$ and
$$ \forall t\in \tilde{U}_{2}\ |\lll  y',T_tx_1\rrr|<1. \eqno(6)$$

Choose a number $n_2\geq n_1$ such that $\tilde{U}_{2}\subset
[0,n_2]$.

Put $x_2=x_1\pm\frac{ y_{n_2}}{10}$, $x_2'=x_1'\pm\frac{
y'_{n_2}}{10}$. We will decide later which pair of signs $ \pm $ to
choose from the four possible cases. Now, we show that
$$ \forall t\in U_1 \ |\lll x_2', T_tx_2  \rrr|>\frac{9}{10}-\frac{3}{10}. \eqno(7)$$

We have $\lll x_2', T_tx_2  \rrr=\lll x_1', T_tx_1 \rrr\pm \lll
x_1', \frac{T_t y_{n_2}}{10} \rrr\pm\lll \frac{ y_{n_2}'}{10},T_tx_2
\rrr $ for any $t$. Therefore,
$$ |\lll x_2', T_tx_2 \rrr|\geq|\lll x_1', T_tx_1
\rrr|-\frac{S_1(t)}{10},\ \ S_1(t)=\left(|\lll x_1', T_t y_{n_2}
\rrr|+|\lll  y_{n_2}',T_tx_2 \rrr|\right) .\eqno(8)$$

If $t\in U_1$ then $S_1(t)<3$. Indeed, if $t\in U_1$ then $|\lll
x_1', T_t y_{n_2} \rrr|\leq\|T_t y_{n_2}\|\leq\frac{11}{10}$ and
$|\lll y_{n_2}',T_tx_2 \rrr|\leq\|T_tx_2\|=\|T_t( y_1\pm\frac{
y_{n_2}}{10})\|\leq \frac{11}{10}+\frac{11}{10^2}$. Recalling (5)
yields (7).

If $t\notin U_{1}$ then (8) is useless for the estimation of $ |
\lll x_2 ', T_tx_2\rrr |$ from below because, for example, we cannot
estimate the value of $|\lll x_1',T_tx_1\rrr|$ from below. Let's
utilize another trick will be called "choosing from the four".

Note that
 $2\frac{ y'_{n_2}}{10}=(x_1'+\frac{ y'_{n_2}}{10})-(x_1'-\frac{ y'_{n_2}}{10})$
and, analogously, $2\frac{ y_{n_2}}{10}=(x_1+\frac{
y_{n_2}}{10})-(x_1-\frac{ y_{n_2}}{10}).$ Therefore,
$$\forall t\ \
4\left|\lll\frac{ y'_{n_2}}{10},T_t\frac{ y_{n_2}}{10}
\rrr\right|\leq\sum_{\pm\in\{+,-\}}\left|\lll (x_1'\pm\frac{
y'_{n_2}}{10}),T_t(x_1\pm\frac{ y_{n_2}}{10})\rrr\right|,
$$ and for each $t\in \tilde{U}_{2}$ at least one of the four terms on the right is at
least $\frac{1}{10^2}|\lll y'_{n_2}, T_t y_{n_2}\rrr|$. At the same
time, by (4), for all $t\in \tilde{U}_2\subset [0,n_2]$\ $|\lll
y'_{n_2}, T_t y_{n_2}\rrr|>\frac{9}{10}$. Therefore, we can choose a
subset $U_2\subset \tilde{U}_{2}$ whose measure $\mu(U_2)\geq
\frac14\mu(\tilde{U}_{2})\geq m_2$ and for some pair of signs $\pm$
(say, for ``$++$'') we have
$$ \forall t\in U_2 \ \left|\lll(x_1'+\frac{ y'_{n_2}}{10}),T_t(x_1+\frac{ y_{n_2}}{10})\rrr\right|=
|\lll x_2', T_tx_2 \rrr|>\frac{9}{10^3}.\eqno (9) $$

{\bf The construction of the vector $x_3$}. Let  $\tilde{U}_3$ be a
compact set, $\mu(\tilde{U}_{3})\geq 4m_3$, $\forall t\in
\tilde{U}_{3}$ $|\lll y', T_t x_2 \rrr|<1.$ Return for a while to
the set $U_2$.  From (6), the compactness of the set $ \{T_tx\mid
t\in \tilde {U} _ {2} \} \subset X $ and the fact that $ y_n
'\overset {\sigma ^ *} {\to} y ' $   it follows that

$$ \exists n_3\geq n_2 \mid \tilde{U}_3\subset [0,n_3],\ \forall n\geq n_3,
 \forall t\in U_2\ |\lll  y'_n,T_tx_1\rrr|<1. \eqno(10)$$

Put $x_3=x_2\pm\frac{ y_{n_3}}{10^3}$, $x_3'=x_2'\pm\frac{
y'_{n_3}}{10^3}$. The pair of signs $ \pm $ will be chosen later.
Now, arguing in the same way as in deriving (8), we obtain:
$$\forall t\ |\lll x_3',T_tx_3\rrr|\geq |\lll
x_2',T_tx_2\rrr|-\frac{S_2(t)}{10^3},\ S_2(t)= |\lll x_2',T_t
y_{n_3}\rrr|+ |\lll y'_{n_3}, T_t x_3\rrr|.\eqno(11)
$$
Show that $\forall t\in U_1\cup U_2$ $S_2(t)<3$. Clearly, $ \|x_2 '
\|\leq \frac{11}{10}$. Arguing as in estimating $S_1 (t) $, we have:
$\forall t\in U_1\
S_2(t)\leq\left(\frac{11}{10}\right)^2+\frac{11}{10}+\frac{11}{10^2}+\frac{11}{10^4}<3.$

If $t \in U_2$ then the argument is a bit more difficult. The first
summand in $S_2 (t) $ is estimated in the old way: if $t\in
U_2\subset [0, n_3] $, then
 $\|T_t y_{n_3}\|\leq
\frac{11}{10}$ and $|\lll x_2',T_t y_{n_3}\rrr|\leq
(\frac{11}{10})^2$. Consider the summand $|\lll y'_{n_3}, T_t
x_3\rrr|$. Recall that $T_tx_3=(T_tx_1+\frac{T_t
y_{n_2}}{10}\pm\frac{T_t y_{n_3}}{10^3}).$

The value $\|\frac{T_t y_{n_2}}{10}\pm\frac{T_t y_{n_3}}{10^3}\|$
with $t\in U_2$ is less than $\frac{11}{10^2}+\frac{11}{10^4}$.

The vector $T_tx_1=T_ty _ {n_1} $, which for $t\in U_2$ can a priori
become large in norm, could have ruined everything  but, owing to
(10), the value of the $ y ' _ {n_3} $ at this vector at $t\in U_2$
is less than $1$. Therefore, $ \forall t\in U_2 \ \left|\lll
y'_{n_3}, T_t x_3\rrr\right| \leq
\left(1+\frac{11}{10^2}+\frac{11}{10^4}\right)$. Finally, $\forall
t\in U_2 \ S_2(t)<
\left(\frac{11}{10}\right)^2+\left(1+\frac{11}{10^2}+\frac{11}{10^4}\right)<3.
$

Now we conclude from (11), (7), (9) that
$$\forall t\in U_1  \ |\lll x_3', T_tx_3  \rrr|>\frac{9}{10}-\frac{3}{10}-\frac{3}{10^3},
\ \ \ \ \forall t\in U_2 \  |\lll x_3', T_tx_3
\rrr|>\frac{9}{10^3}-\frac{3}{10^3}.  $$ \vskip1mm

Choose of the pair ``$\pm$'' and the set $U_3\subset\tilde{U}_3$,
$\mu(U_3)\geq m_3$ with the help of the  ``choosing from the four''
trick. We have:
$$ \forall t\in U_3 \ |\lll x_3', T_tx_3 \rrr|>\frac{9}{10}\cdot\left(\frac{1}{10^3}\right)^2=\frac{9}{10^7}. $$

Note that it was quite easy to construct the vector $x_1$. At the
step 2 we need the ``choosing from the four'' trick. The novelty of
the step 3 was the usage of (10) and the preparation for it --- the
formula (6) --- should be made at the very beginning of the step 2.
Next steps have no significant differences from the step 3.

{\bf The construction of the vector \it $x_l$}\rm, $l\geq 3$.
Suppose that we have constructed the sets $U_1, \ldots, U_l $, the
numbers $n_1\leq n_2\leq\cdots\leq n_l $, the vectors $x_1, \ldots,
x_l $ of the form $x_l =\sum _ {i=1} ^l \pm\frac {y _ {n_i}} {10 ^
{2 ^ {i-1}-1}} $ and, similarly, $x_1 ', \ldots, x_l ' $ so that the
following properties hold:

$\begin{array}{llr} 1_l) & U_i\subset [0,n_i], \mu(U_i)\geq m_i,
& i=1,2,\ldots,l;\\
2_l) & \forall t\in U_i\  |\lll  y', T_tx_{i-1}\rrr|<1,&
i=2,\ldots,l;\\
3_l) & \forall t\in U_i  \ \forall n\geq n_{i+1}\ |\lll  y_n',
T_tx_{i-1}\rrr|<1, & i=2,\ldots, l-1;\\
4_l) & \forall t\in U_i\ |\lll x_l',T_tx_l\rrr|\geq
\frac{9}{10^{2^i-1}}-3\sum_{j=i}^{l-1} \frac{1}{10^{2^j-1}}, &
i=1,\ldots ,l-1;\\
5_l) & \forall t\in U_l\ |\lll x_l',T_tx_l\rrr|\geq
\frac{9}{10^{2^l-1}}.\\
\end{array}$

Now we  construct a set $U _ {l+1} $, a number $n _ {l+1} $, a
vector $x _ {l+1} $ and the corresponding $x _ {l+1} ' $ so that the
properties $1_{l+1})-5_{l+1})$ also hold.

Choose a compact set $ \tilde {U} _ {l+1} $ such that $ \mu (\tilde
{U} _ {l+1}) \geq 4m _ {l+1} $, $ \forall t\in \tilde {U} _ {l+1}\ |
\lll y ', T_t x_l\rrr | <1$. The possibility of such choice follows
from  (3).

Choose $n_{l+1}\geq n_l$  such that
$$\tilde{U}_{l+1}\subset[0,n_{l+1}],\ \forall n\geq n_{l+1}\ \forall
t\in U_l\ |\lll y_n', T_tx_{l-1}\rrr|<1.$$ Such $n _ {l+1} $ exists
by condition $2_l $) (see the argument before~(10)). It is condition
$3 _ {l+1} $.

Put $x_{l+1}=x_l\pm\frac{ y_{n_{l+1}}}{10^{2^l-1}}$,
$x_{l+1}'=x_l'\pm\frac{ y_{n_{l+1}}'}{10^{2^l-1}}$.  Choose a pair
of signs
 $\pm$ and a subset $U_{l+1}\subset\tilde{U}_{l+1}$ using the ``choosing from the four''
trick.  Thus, conditions $1_{l+1}$, $2_{l+1}$ hold as well as
$5_{l+1}$:

 $$\forall t\in U_{l+1}\ |\lll x_{l+1}', T_tx_{l+1}\rrr|\geq \left|\lll\frac{ y_{n_{l+1}}'}{10^{{2^l-1}}},
 \frac{T_t y_{n_{l+1}}}{10^{{2^l-1}}}\rrr\right|\geq\frac{9}{10}\cdot\left(\frac{1}{10^{{2^l-1}}}\right)^2=
 \frac{9}{{10^{{2^{l+1}-1}}}}
$$

The last, check the condition $4_{l+1}$. For all $t$
$$
|\lll x_{l+1}', T_tx_{l+1} \rrr| \geq |\lll x_{l}', T_tx_{l} \rrr|-
\frac{S_l(t)}{10^{2^l-1}},\ \
 S_l(t)=|\lll x_l',T_t y_{n_{l+1}}
\rrr|+|\lll  y_{n_{l+1}}', T_t x_{l+1}\rrr|.$$

It suffices to establish that $S_l (t) <3$ for all $t\in
U_1\cup\cdots\cup U_l$. The first summand $S_l (t) $ is less than
$\|x_l'\|\cdot\|T_ty_{n_{l+1}}\|\leq\frac {12\cdot11} {100} $. Let
us estimate the second summand.

If $t\in U_i\subset [0,n_i]$, then, writing
$x_{l+1}=x_{i-1}+\sum_{j=i}^{l+1}(\pm\frac{
y_{n_j}}{10^{2^{j-1}-1}})$, we have

$$ |\lll  y'_{n_{l+1}}, T_t x_{l+1}\rrr|=\left|\left\lll  y'_{n_{l+1}}, T_tx_{i-1}+\sum_{j=i}^{l+1}(\pm
\frac{T_t y_{n_j}}{10^{2^{j-1}-1}}) \right\rrr\right| < |\lll
y'_{n_{l+1}}, T_tx_{i-1}\rrr|+\frac12<\frac32.
$$

The inequalities in the previous formula are valid due to the
smallness of $ \|T_t y_j \| $ for $j\geq i $ (as $t\in
U_i\subset[0,n_i]$), and also the condition $3 _ {l+1}, $ which has
been already proved above (sf. the estimation of $S_2(t)$ with $t\in
U_2$.)

Let $x'=\lim x_k'$ and $x=\lim x_k$. Then,
 by "$4_\ii$",  we have
$$\forall t\in U_i\ |\lll x',T_tx\rrr|\geq
\frac{9}{10^{2^i-1}}-3\sum_{j=i}^{\infty}
\frac{1}{10^{2^j-1}}>\frac{5}{10^{2^i-1}} .$$

The theorem is proved.

\section{On the
possibility of choosing a smooth vector $x$ in Theorem 1}

Is it possible to choose the vector $x$ to be smooth in Theorem 1?
Clearly, it cannot be done if, for example, $\omega_1<s_0=0$: in
this case the semigroup is ES; therefore, the orbits of smooth
vectors decrease too rapidly. Based on example from \cite{17Zb},
Wrobel in \cite{18Wr} constructed a semigroup with
$s<\omega_1<s_0=\omega_0$, moreover, in his example
$\omega_n=2^{-n}$, where $\omega_n$ is the growth of the semigroup
on $D(A^n)$. So, this semigroup, after having been rescaled so that
$\omega_1=0$, remains unbounded, sytisfies the Theorem 1 but
decreases exponentially on any smooth vector (and the more
smoothness, the faster decrease). The semigroup from \cite{17Zb}
itself, renormed so that   $s=-1$, $\omega_1=s_0=\omega_0=0$, is not
even ES; however, (\cite{19BEKP}): $\|T_t x\|=O(1/t)$ for all $x\in
D(A^{1+\ee})$ $\forall \ee>0$.

 What can impede a vector $x=\sum \gamma_j y_{n_j}$ of Theorem 1
from being in $D(A)$ if the number series $\sum\gamma_i\in\Bbb R$
converges and $y_{n_j}\in D^{\ii}(A)$, $\|y_{n_j}\|=1$? The answer
is simple: the series $\sum \gamma_j Ay_{n_j}$ does not have to
converge in $X$, as the norms $\|Ay_{n_j}\|\approx|\beta_{n_j}|$ can
grow fast with $j$ (see the proof of the Lemma 1).

Suppose that $s_0=s=0$ and the bound of Sp(A) is reached, i.e. there
exists $\lambda\in Sp(A),\ Re\lambda=0$ (a typical case for
$C_0$-semigroups). Then we may take the numbers $i\beta_n$ in Lemma
1 near this $\lambda$ (not somewhere "near infinity") and, the norms
$\|Ay_{n_j}\|\approx|\beta_{n_j}|$ will be bounded with $j$. Then
the series $\sum \gamma_j Ay_{n_j}$ converges. Operator $A$ is
closed, so $x$  gets to $D(A)$. Let us prove that it is even
possible to find \it infinitely \rm smooth vector $x$:
\smallskip

{\bf Theorem  2}\it. Suppose that $s(A)\geq 0$ is reached. Then
there exist $x'\in X'$, $x\in D(A^\ii)$ satisfying the hypothesis of
Theorem 1.

\it{Proof}\rm. Let $\lambda\in Sp(A),\ Re(\lambda)=s$. After
rescaling by $e^{-\lambda}$, we may assume $s=0\in Sp(A)$. It
suffices to find the elements $y_n$ of the proof of Theorem 1 such
that $y_n\in D(A^\ii)$ and $\sum_n \gamma_n A^k y_n$ converges in
$X$ for all $k=0,1,2,\ldots$. In the next Lemma we will show that it
can be done. This Lemma is of interest on its own.\smallskip

\bf Lemma 2\it. If  $0\in Sp(A)$, then $\forall \delta>0$ $\forall
n\in\Bbb N$ $\exists y_n\in D(A^\ii)$ such that
$$\|y_n\|=1, \ \forall i\leq n \ \ \|A^iy_n\|<\delta.\eqno(12_n)$$

Proof\rm.  Let $n=1$.  Operator $(\lambda-A)^{-1}$ is unbounded with
$\lambda$ near $0$, so, $12_1$) followed by  Lemma 1.

For $n>1$ we use the construction of the scale of associated Sobolev
semigroups (\cite{14Na}, A-I 3.5). On the space $D(A^n)$, consider
the ``iterated'' graph norm $\|x\|_n=\|x\|+\|Ax\|+\cdots+\|A^nx\|$.
The semigroup $T_t:D(A^{n-1})\to D(A^{n-1})$ with the generator
$A:D(A^{n})\to D(A^{n-1})$ is isomorphic to the initial one.

Rewriting ($12_1$) for this semigroup, we infer that there exists
$y_n\in D(A^\ii)$ such that
$$ \|y_n\|+\|Ay_n\|+\cdots+\|A^{n-1}y_n\|=1,\
\|Ay_n\|+\cdots+\|A^{n}y_n\|<\delta.\eqno(13)$$

It follows from the second expression (inequality) of (13) that
$\|Ay_n\|<\delta,$ $\|A^2y_n\|<\delta$,\ldots,
$\|A^{n}y_n\|<\delta$. But then from the first \it equality \rm of
(13) we obtain that $\|y_n\|$ is about 1. The rest is obvious. Lemma
3 and Theorem 2 are proved.\smallskip

Note that Theorem 2 seems quite natural, if we take into account the
semigroup of shifts: clearly, a function on $[0,\ii)$  can be made
to decrease arbitrary slowly and the infinite differentiability, as
a local phenomenon, is not an obstacle here. \vskip1mm

\bf Acknowledgments\rm.  We express our gratitude to Professor Jan
van Neerven for his interest in our work and stimulating questions.
In particular, his advice enabled us to substantially strengthen the
main results.

We thank Professor Constantin Buse for information that stimulated
our interest in the problem.

Our thanks also to Professor Bogdan Sasu who sent me a copy of a
rather inaccessible paper \cite{16Fo}.

\end{document}